\documentclass{article}

\usepackage{arxiv}

\usepackage[utf8]{inputenc} % allow utf-8 input
\usepackage[T1]{fontenc}    % use 8-bit T1 fonts
\usepackage{url}            % simple URL typesetting
\usepackage{booktabs}       % professional-quality tables
\usepackage{amsfonts}       % blackboard math symbols
\usepackage{nicefrac}       % compact symbols for 1/2, etc.
\usepackage{microtype}      % microtypography
\usepackage{lipsum}         % Can be removed after putting your text content
\usepackage{graphicx}
\usepackage{subcaption}
\usepackage{doi}

\usepackage{colortbl}
\usepackage{listings}
\usepackage{bm}
\usepackage{physics}

\usepackage{amsmath,accents}
\usepackage{cleveref}       % smart cross-referencing

\usepackage[ruled,vlined]{algorithm2e}

\usepackage{appendix}
\usepackage{cite}

\usepackage{multirow}
\usepackage{array}

\usepackage{chngcntr}

\usepackage{dcolumn}

\title{Physics-informed neural networks and neural operators for a study of EUV electromagnetic wave diffraction from a lithography mask}

\date{}

\author{ \href{https://orcid.org/0000-0002-4930-1846}{\includegraphics[scale=0.06]{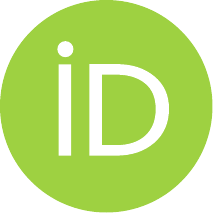}\hspace{1mm}Vasiliy A. Es'kin}\thanks{Corresponding author: Vasiliy Alekseevich Es’kin (\href{vasiliy.eskin@gmail.com}{vasiliy.eskin@gmail.com})} \\
	Department of Radiophysics, University of Nizhny Novgorod\\
	23 Gagarin Ave., Nizhny Novgorod 603022, Russia\\
	\href{vasiliy.eskin@gmail.com}{\texttt{vasiliy.eskin@gmail.com}} \\
	\And
	{Egor V. Ivanov} \\
	Department of Radiophysics, University of Nizhny Novgorod\\
	23 Gagarin Ave., Nizhny Novgorod 603022, Russia\\	\texttt{iev90078@gmail.com}
}

\renewcommand{\headeright}{}
\renewcommand{\undertitle}{}
\renewcommand{\shorttitle}{A preprint}

\renewcommand{\vec}{\bf}
\newcommand{\eps}{\varepsilon}

\begin {document}

\maketitle

\index{Es'kin, V.A.}
\index{Ivanov, E.V.}

\begin{abstract}	
	Physics-informed neural networks (PINNs) and neural operators (NOs) for solving the problem of diffraction of Extreme Ultraviolet (EUV) electromagnetic waves from a mask are presented. A novel hybrid Waveguide Neural Operator (WGNO) is introduced, which is based on a waveguide method with its most computationally expensive part replaced by a neural network. Numerical experiments on realistic 2D and 3D masks show that the WGNO achieves state-of-the-art accuracy and inference time, providing a highly efficient solution for accelerating the design workflows of lithography masks.
\end{abstract}

\section{Introduction}

The relentless advancement of semiconductor technology, as described by Moore's law, has led the feature sizes of integrated circuits into the nanometer scale. Extreme Ultraviolet lithography, utilizing wavelengths around 13.5 nm, is the cornerstone technology for manufacturing these state-of-the-art semiconductor devices. At such scales, the diffraction of electromagnetic waves from the photomask becomes a dominant physical effect, causing the pattern printed on the silicon wafer to deviate significantly from pattern of the mask.

In processor manufacturing, Optical Proximity Correction (OPC) techniques are essential. OPC involves intentionally distorting the mask pattern to counteract the anticipated diffraction effects, thereby producing the desired circuit pattern on the wafer. A critical step in the OPC loop is the forward simulation: accurately and rapidly predicting the electromagnetic field distribution after diffraction from a complex multi-layered EUV mask.

Traditional methods for this simulation task present a trade-off. Simple models like the infinitely thin opaque mask approximation are fast but suffer from low accuracy, which is unacceptable for modern manufacturing. Rigorous numerical solvers such as the Finite Element Method (FEM) or the Waveguide (WG) method provide high accuracy but are computationally expensive. The WG method is among the fastest and most accurate rigorous methods, yet it can still be a bottleneck in the design process~\cite{it1}. Recently, deep neural networks, such as Convolutional Neural Networks (CNNs)~\cite{it1} or U-Net~\cite{it2}, have been proposed, providing extremely fast inference. However, neural network approaches are often based on supervised learning, require a large dataset, have a significant training time, and often do not demonstrate the necessary accuracy of the solution and the degree of its generalization.

This work investigates an alternative deep learning paradigm for the given problem. We apply physics-informed neural networks \cite{it3,it4} and propose a novel neural operator to solve the EUV diffraction problem. These neural networks are trained in an unsupervised manner, leveraging the governing physical equations directly in the training process. We evaluate their performance on benchmark problems and on realistic 2D and 3D mask models for current and promising industrial lithography systems (13.5 and 11.2 nm) \cite{it6,it7}.

\section{Problem Formulation}

\begin{figure}[ht!]\centering
	\includegraphics[width=.43\textwidth]{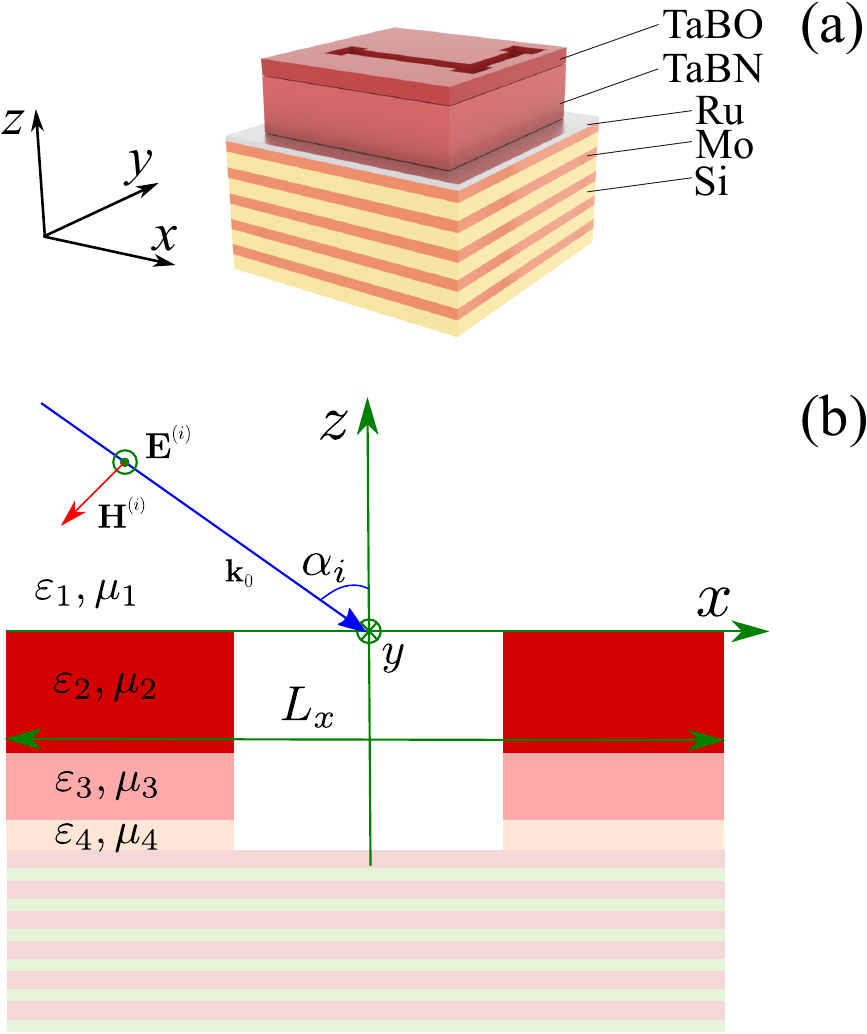}
	\caption{Geometry of the problem in the stereometric view (a) and the cross section $y$ = 0 (b).}\label{fig1}
\end{figure}

Consider a layered mask of thickness $D$, consisting of $J$ layers, located in free space and enclosed in an interval $[-D,0]$ along the $z$-axis of a Cartesian coordinate system ($x,y,z$) , as shown in Fig.~\ref{fig1}. The mask structure is periodic along the $x$ and $y$ axes with periods $L_x$ and $L_y$, respectively. Each layer of the mask is filled with a medium, which is uniform in the $z$ direction and has a dielectric permittivity $\eps_j(x,y)$. The dielectric permittivities of the media were obtained from experimental data published in~\cite{HENKE1993181,CenterXRayOpt}.

We consider the diffraction of an electromagnetic monochromatic plane wave with angular frequency $\omega$ from a mask (see Fig.~\ref{fig1}(b)). The  magnetic fields in the incident wave are given, with $\exp(i\omega t)$ time dependence dropped, by ${\vec H}^{(i)}={\vec H}_0\exp[-i \left(k_{0;x} x + k_{0;y} y - k_{0;z} z\right)]$, where ${\vec H}_0$ is the magnetic field amplitude, $k_{0;x}$, $k_{0;y}$, and $k_{0;z}$ are the components of the wave vector ${\vec k}_0$ in free space ($k_0 = \left(k_{0;x}^2 + k_{0;y} ^ 2 + k_{0;z}^2\right)^{1/2}$, $k_0 = \omega /c$, where $c$ is the speed of light in free space), and the superscript $(i)$ denotes the incident wave. To simplify further consideration of the problem, we assume that the field of incident wave is periodic along the $x$ and $y$ axes with periods $L_x$ and $L_y$, respectively.

We assume that the media of the mask are time-independent and nonmagnetic. It can be shown based on Maxwell equations that the vector components $\vec{E}$ and $\vec{H}$ of electromagnetic field in a homogeneous medium along the $z$ axis within each $j$th layer are expressed in terms of magnetic components $H_x$ and $H_y$ which, in turn, satisfy the following system of equations:
\begin{align}
	&\Delta H_x + k_0^2 \varepsilon_j H_x + \frac{1}{\varepsilon_j} \frac{\partial \varepsilon_j}{\partial y} \left( \frac{\partial H_{y}}{\partial x} - \frac{\partial H_{x}}{\partial y} \right){=} 0,\notag \\
	&\Delta H_y + k_0^2 \varepsilon_j H_y +  \frac{1}{\varepsilon_j}\frac{\partial \varepsilon_j}{\partial x} \left( \frac{\partial H_{x}}{\partial y} - \frac{\partial H_{y}}{\partial x} \right){=} 0. \label{eq1}
\end{align}
The problem is solved within a computational domain subject to boundary conditions. Periodicity conditions are satisfied on the lateral boundaries (${\vec{H}}\left(-L_x/2, z\right) = {\vec{H}}\left(L_x/2, z\right)$, ${\vec{H}}\left(-L_y/2, z\right) = {\vec{H}}\left(L_y/2, z\right)$). At the top ($z=z_{\max}$) and bottom ($z=z_{\min}$) of the domain, Sommerfeld radiation conditions are applied to ensure that scattered waves propagate outwards without unphysical reflections. These are formulated for every spatial Fourier harmonic of the field as follow:
\begin{equation}\label{eq:sommerfeld}
	\int\limits_{0}^{L_x} \int\limits_{0}^{L_y} \left(\frac{\partial {{\vec H}_\perp}}{\partial z} \pm i {k}_{z;mn} {{\vec H}_\perp} \right) e^{i (\kappa_x m x + \kappa_y n y)} dx dy = 0,
\end{equation}
where ``$+$'' and ``$-$'' correspond to $z=z_{\max}$ and  $z=z_{\min}$, respectively, ${k}_{z;mn} = (k_0^2 - \kappa_x^2 m^2 - \kappa_y^2 n^2)^{1/2}$ (branch ${\rm Re}\,({k}_{z;mn}) \geq 0$ is taken), $\kappa_x = 2 \pi / L_x$, $\kappa_y = 2 \pi / L_y$, $m$ and $n \in \mathbb{Z} $.

\section{Solution Methods}

\subsection{Finite Element Method (FEM)}
The Finite Element Method is a powerful and general numerical technique for solving partial differential equations. The core principle of FEM involves discretizing the computational domain into a mesh of smaller, simpler subdomains called finite elements. Within each element, the unknown physical field is approximated by a set of basis functions. Assembling the equations for all elements leads to a large global system of linear equations, which is then solved to obtain the approximate solution across the entire domain. While robust, FEM can be computationally expensive for problems like EUV diffraction, which require very fine meshes to resolve the small wavelength features. In our comparative studies, we utilize the open-source solver FreeFem++~\cite{hecht2012new}.

\subsection{Waveguide Method}

The Waveguide (WG) method serves as our high-fidelity reference solver~\cite{it1,lucas1996efficient}. It leverages the layered structure of the mask. Within each layer, the field is represented as a superposition of waveguide modes. The solution is sought in the form $H_{x,y}(x,y,z) = h_{x,y}(x,y) \exp(-ik_z z)$. By expanding the transverse functions $h_{x,y}$ and the permittivity $\varepsilon(x,y)$ into Fourier series, the governing PDEs~(\ref{eq1}) are transformed into a large algebraic eigenvalue problem for each layer $j$ (see~\cite{lucas1996efficient}):
\begin{align}
	\mathbf{M}^{(j)}_{\text{layer}} \begin{pmatrix} \mathbf{B}^{(j)} \\ \mathbf{C}^{(j)} \end{pmatrix} = \left(k_z^{(j)}\right)^2 \begin{pmatrix} \mathbf{B}^{(j)} \\ \mathbf{C}^{(j)} \end{pmatrix},
\end{align}
where $\mathbf{M}^{(j)}_{\text{layer}}$ is the matrix derived from PDEs~(\ref{eq1}), $\mathbf{B}^{(j)}$ and $\mathbf{C}^{(j)}$ are the eigenvectors which are vectors of Fourier coefficients for $h_x$ and $h_y$. The eigenvalues $\left(k_z^{(j)}\right)^2$ give the propagation constants. The total field in each $j$th layer is a linear combination of modes
\begin{align}
	\begin{bmatrix} H_x^{(j)} \\ H_y^{(j)} \end{bmatrix} &{=} \sum\limits_{p = 1}^{2N} k_{z;p}^{(j)} \left[A_{p;1}^{(j)} e^{i k_{z;p}^{(j)} z} + A_{p;2}^{(j)} e^{-i k_{z;p}^{(j)} z}\right] \sum\limits_{m,n{=}-N_x,-N_y}^{N_x,N_y} \begin{bmatrix}{B}_{p,mn}^{(j)} \\ {C}_{p,mn}^{(j)}\end{bmatrix} \psi_{mn},
\end{align}
where $\psi_{mn} = \exp\left(- i \kappa_x m x - i \kappa_y n y\right)$, $N_x$ and $N_y$ are maximum numbers of the Fourier harmonics along the $x$ and $y$ axes, $N = (2 N_x + 1) (2 N_y + 1)$. The reflected field ($z>0$) is written as:
\vspace{-0.1cm}
\begin{align}
	\begin{bmatrix} H_x^{(r)} \\ H_y^{(r)} \end{bmatrix} {=} \hspace{-1mm}\sum\limits_{m,n = -N_x,-N_y}^{N_x,N_y} \hspace{-3mm}{k}_{z;mn} \begin{bmatrix} A_{x;mn}^{(r)} \\ A_{y;mn}^{(r)} \end{bmatrix}  \psi_{mn} e^{- i {k}_{z;mn} z}.
\end{align}
The transmitted field ($z<-D$) is written as:
\vspace{-0.1cm}
\begin{align}
	\begin{bmatrix} H_x^{(t)} \\ H_y^{(t)} \end{bmatrix} {=} \hspace{-1mm}\sum\limits_{m,n = -N_x,-N_y}^{N_x,N_y} \hspace{-3mm}{k}_{z;mn} \begin{bmatrix} A_{x;mn}^{(t)} \\ A_{y;mn}^{(t)} \end{bmatrix}  \psi_{mn} e^{i {k}_{z;mn} z}.
\end{align}
Here $A_{p;1}^{(j)}$, $A_{p;2}^{(j)}$, $A_{x;mn}^{(r)}$, $A_{y;mn}^{(r)}$, $A_{x;mn}^{(t)}$ and $A_{y;mn}^{(t)}$ are the unknown coefficients. By satisfying continuity conditions of the tangential field components at each layer interface, a global system of linear equations is formed:
\begin{equation} \label{eq:wg_linear_system}
	\hat{\mathbf{M}} \mathbf{A} = \mathbf{R},
\end{equation}
where $\hat{\bf M}$ is matrix of system of equations obtained from the boundary conditions, $\mathbf{A}$ is the vector of unknown mode amplitudes (noted above) and $\mathbf{R}$ is determined by the incident field. Solving this large linear system is the most computationally expensive part of the WG method.

\subsection{Physics-Informed Neural Networks}
A PINN approximates the solution of problem, e.g., $H_x(\mathbf{x})$ and $H_y(\mathbf{x})$, with deep neural networks $H_x(\mathbf{x};\bm{\theta})$ and $H_y(\mathbf{x};\bm{\theta})$, where $\mathbf{x}=(x,y,z)$ are the spatial coordinates and $\bm{\theta}$ are the trainable network parameters (weights and biases). The network is trained by minimizing a composite loss function that enforces the physical laws:
\begin{equation}
	\mathcal{L}(\bm{\theta}) = \lambda_{bc} \mathcal{L}_{bc} + \lambda_{r} \mathcal{L}_{r}.
\end{equation}
$\mathcal{L}_{bc}$ is a boundary loss term that penalizes deviations from the boundary conditions (\ref{eq:sommerfeld}), and $\mathcal{L}_{r}$ is a residual loss term that penalizes non-zero residuals of the governing PDEs~(\ref{eq1}) (see details in~\cite{it3,it4}), $\lambda_{bc}$ and $\lambda_{r}$ are hyperparameters, which allow for separate tuning of {the} learning rate for each of the loss terms in order to improve the convergence of the model (see the procedure for choosing these parameters and examples~\cite{it3,it4,Eskin2025}). The derivatives required to compute the PDE residuals are obtained via automatic differentiation. The optimization problem can be defined as follows:
\begin{equation}
	{\bm \theta}^* = {\arg}\,\underset{{\bm \theta}}{\min} \, \mathcal{L}({\bm \theta}),\label{eq9}
\end{equation}
where ${\bm \theta}^*$ are optimal parameters of the neural network which minimize the discrepancy between the exact unknown solution $H_x(\mathbf{x})$ and $H_y(\mathbf{x})$ and the approximate one $H_x(\mathbf{x};\bm{\theta}^*)$ and $H_y(\mathbf{x};\bm{\theta}^*)$. This approach does not require any pre-computed data for training.

\subsection{Waveguide Neural Operator}
To overcome the limitations of both the WG method (slow solving the linear system) and PINNs (low accuracy for complex geometries), we propose a hybrid Waveguide neural operator (WGNO). The core idea is to replace only the most computationally intensive part of the WG method with a neural network. Specifically, we train a multi-layer perceptron (MLP) to learn the mapping from the right-hand side vector $\mathbf{R}$ and the mask parameters (Fourier coefficients of permittivity) to the solution vector $\mathbf{A}$ of the linear system (\ref{eq:wg_linear_system}):
\begin{equation}
	\mathbf{A}_{\bm{\theta}} = \text{MLP}\left(\mathbf{R}, \{\varepsilon_{mn}^{(j)}\}^{J}_{j=1}; \bm{\theta}\right),
\end{equation}
where
\begin{align}
	& {\eps}^{(j)}_{mn} = \frac{1}{L_x L_y} \int\limits_{0}^{L_x} \int\limits_{0}^{L_y} \eps_j(x,y) e^{i \kappa_x m x + i \kappa_y n y} dx dy,
\end{align}
$\{\varepsilon_{mn}^{(j)}\}^{J}_{j=1}$ is set of Fourier coefficients of permittivities. The rest of the WG pipeline, such as calculating the modes, forming $\mathbf{R}$, and reconstructing the field in real space from the predicted $\mathbf{A}_{\bm{\theta}}$, remains unchanged. The network is trained to minimize the following loss function: $\mathcal{L}(\bm{\theta}) = \|\hat{\mathbf{M}} \mathbf{A}_{\bm{\theta}} - \mathbf{R}\|^2_2$. This hybrid approach operates in the mesh-independent Fourier (latent) space, inherits the physical structure of the WG method, and directly targets the primary computational bottleneck.

\section{Numerical Experiments}

We evaluated the performance of the proposed methods on several test problems. The neural networks were implemented in PyTorch. For PINNs, we used a MLP with 3 hidden layers and 128 neurons each with hyperbolic tangent activation function ($\tanh$). For WGNO, the MLP had 2 hidden layers.

\subsection{Validation on Test Problems}
To understand the accuracy of solutions to problems provided by numerical solvers (FreeFem++, WG) and solvers based on artificial neural networks, we have tested the methods on three 2D problems with known analytical solutions: (1) plane wave propagation in a homogeneous medium, (2) reflection from a single interface, and (3) reflection from a dielectric layer. 

We assume the media are homogeneous in the $y$ direction. Electromagnetic wave of TE polarization propagates in free space with wave vector ${\vec k}= k_x^{(i)} {\vec x}_0 - k_z^{(i)} {\vec z}_0$ ($k_x^{(i)}$ and $k_z^{(i)}$ are positive values, and $k_x^{(i)} = k_z^{(i)}$). The medium under interface and in the dielectric layer have permittivity $\varepsilon = 4$. The dielectric layer thickness is $\pi/ k_z^{(i)}$. Coordinates of points $(x, z)$ of area in which we considered the solutions satisfy to relations $x \in [-\pi/ k_x^{(i)}, \pi/ k_x^{(i)}]$ and $z \in [-\pi/ k_z^{(i)}, \pi/ k_z^{(i)}]$ ($z \in [-2\pi/ k_z^{(i)}, \pi/ k_z^{(i)}]$ for the dielectric layer).

For the numerical solver we divided every axis in $100$ intervals, as a result we had about $10^4$ collocation points. $N_x$ equals to $10$ for WG and WGNO methods at all test experiments.

For the PINN and WGNO approach we divided the segment along every axis in $100$ intervals. In our experiments for the PINN, we used two-stage learning which consisted of 1000 epochs of optimization with the Adam optimizer with learning rate $10^{-3}$ and 5 epochs of optimization with the LBFGS optimizer. For WGNO we used 1000 epochs of optimization with the Adam optimizer with learning rate $10^{-3}$ and  $10^{-5}$. We have carried out training with randomly initialized neural networks with the Glorot scheme initialization~\cite{Glorot10a}, performing training 7 times with different seeds.

The results are summarized in Table~\ref{tab:test_problems}. The WG method achieves machine precision, as expected. FreeFem++ provides good accuracy but is orders of magnitude less accurate. The PINN achieves reasonable accuracy, but it is less accurate than FEM and struggles with the more complex third problem. The proposed WGNO demonstrates high performance, with errors several orders of magnitude smaller than the PINN, highlighting the effectiveness of the hybrid physics-based approach.

\begin{table*}[h!]
	\centering
	\caption{Performance comparison on 2D test problems. Inference times are $1.2\times 10^{-3}$ s and $1.7\times 10^{-4}$ s for PINN and WGNO, respectively.}
	\label{tab:test_problems}
	\begin{tabular}{l|cc|cc}
		\toprule
		\multirow{2}{*}{Problem} & \multicolumn{2}{c|}{Waveguide} & \multicolumn{2}{c}{FreeFem++} \\
		& Rel. $L_2$ Error & Time (s) & Rel. $L_2$ Error & Time (s) \\
		\midrule
		1. Homogeneous & $2.8\times 10^{-15}$ & $1.0\times 10^{-3}$ & $2.3\times 10^{-3}$ & 3.47 \\
		2. Interface & $2.7\times 10^{-15}$ & $2.0\times 10^{-3}$ & $3.8\times 10^{-3}$ & 3.39 \\
		3. Layer & $6.1\times 10^{-15}$ & $6.0\times 10^{-3}$ & $8.1\times 10^{-3}$ & 5.03 \\
		\midrule \midrule
		\multirow{2}{*}{Problem} & \multicolumn{2}{c|}{PINN} & \multicolumn{2}{c}{WGNO} \\
		& Rel. $L_2$ Error & Training (s) & Rel. $L_2$ Error & Training (s) \\
		\midrule
		1. Homogeneous & $1.7{\times} 10^{-4}  {\pm} 5.3 {\times} 10^{-5}$ & $855$ & $4.0{\times} 10^{-8}{\pm}
		1.6{\times} 10^{-8}$ & $3.7$ \\
		2. Interface & $1.6{\times} 10^{-3} {\pm} 7.9 {\times} 10^{-4}$ & $1532$ & $1.8{\times} 10^{-7} {\pm} 2.8 {\times} 10^{-7}$ & $3.7$ \\
		3. Layer & $5.5{\times} 10^{-2} {\pm} 4.4{\times} 10^{-3}$ & $1221$ & $4.7{\times} 10^{-5}{\pm}
		5.3 {\times} 10^{-5}$ & $11$ \\
		\bottomrule
	\end{tabular}
\end{table*}

\subsection{2D Lithography Mask Simulation}
\begin{figure}[t!]\centering
	\includegraphics[width=0.5\textwidth]{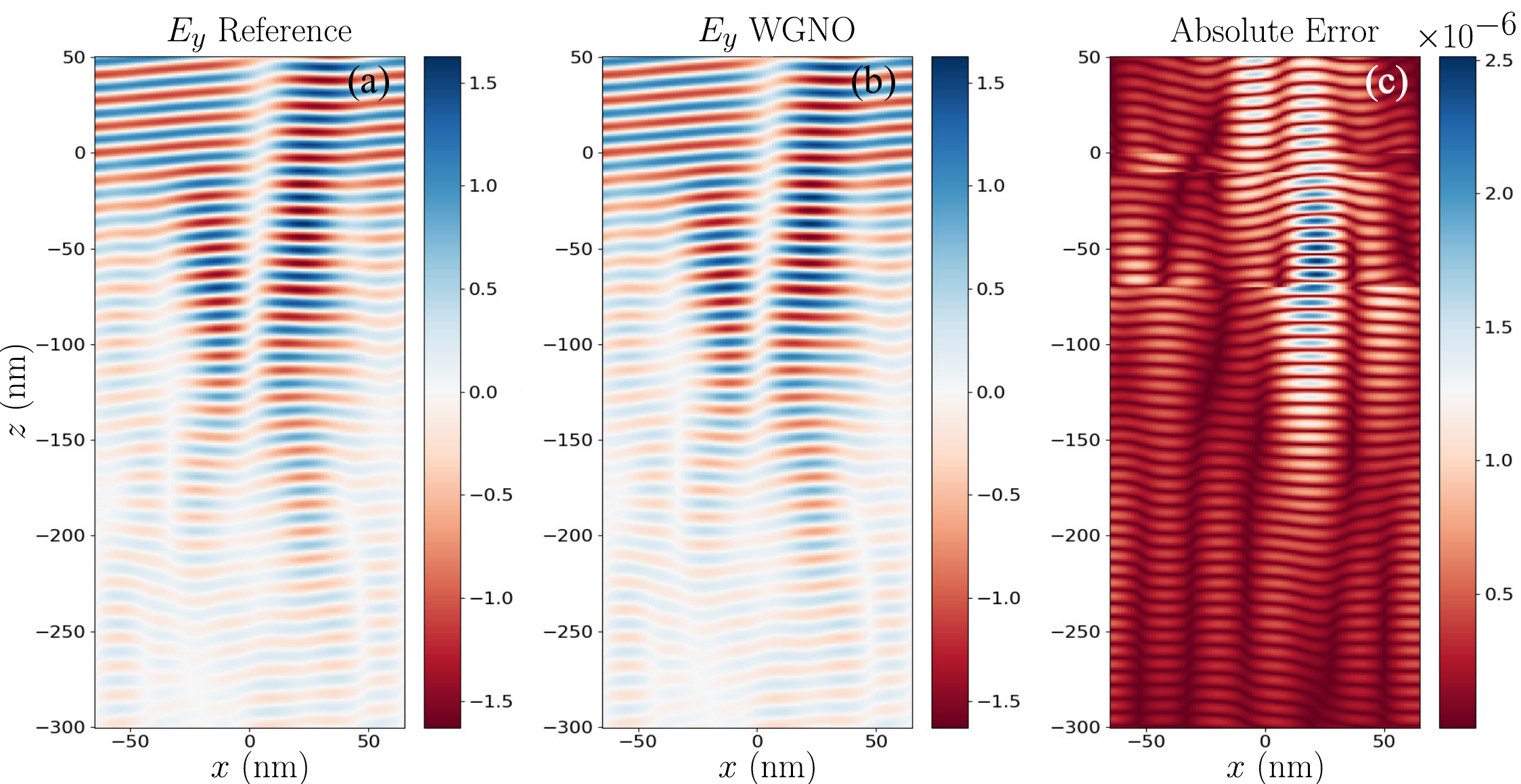}
	\caption{Comparison of reference WG solution (a) and WGNO solution (b) for a 2D mask at 13.5 nm, (c) is absolute error.}\label{fig:2d_mask_results}
\end{figure}

We simulated a realistic 2D EUV mask (see Fig.~\ref{fig1}) for wavelengths ($\lambda$) of 13.5 and 11.2 nm, included an absorber consisting of a TaBO layer (10 nm) and a TaBN layer  (10 nm for PINN, and 60 nm for WGNO), 2.0 nm thin Ru layer,  MoSi layer for PINN and a 31-layer MoSi mirror for the WGNO. The thicknesses of Mo layers were $3$ nm at $\lambda=13.5$ nm and $2.49$ nm at $\lambda=11.2$ nm, of Si layer were $4$ nm and $3.32$ nm for the same wavelengths. The absorbers had a hole with width $L_x/2$. Note, we have obtained the permittivities of media from experimental data~\cite{HENKE1993181,CenterXRayOpt}. The angle of incidence of the wave is $6^{\circ}$. Fig.~\ref{fig:2d_mask_results} shows the comparison of the calculated $E_y$ field, normalized to the amplitude of the incident wave, for the 13.5 nm case. It follows from the calculation results that while the PINN captures the general wave pattern, its absolute error is significant. In contrast, the WGNO result is visually indistinguishable from the reference WG solution, with a very small absolute error (see Fig.~\ref{fig:2d_mask_results}).

The quantitative results are summarized in Table~\ref{tab:2d_mask}. The PINN struggles with this complex problem, yielding errors of several percent after long training times. The WGNO, however, achieves excellent accuracy with a short training time.

\begin{table}[h!]
	\centering
	\caption{Performance on 2D lithography mask simulation. Inference times are $1.2\times 10^{-3}$ s and $1.8\times 10^{-4}$ s for PINN and WGNO, respectively.}
	\label{tab:2d_mask}
	\begin{tabular}{l|ccc}
		\toprule
		$\lambda$ (nm) & Method & \!Rel. \!\!$L_2$\! Error & \! Training (s)\\
		\midrule
		\multirow{2}{*}{13.5} & PINN & $4.9{\times} 10^{-2}$ & $2637$ \\
		& WGNO & $9.5{\times} 10^{-7}$ & $0.012$ \\
		\midrule
		\multirow{2}{*}{11.2} & PINN & $8.9{\times} 10^{-2}$ & $2790$ \\
		& WGNO & $3.9{\times} 10^{-6}$ & $0.017$  \\
		\bottomrule
	\end{tabular}
\end{table}

\subsection{3D Lithography Mask Simulation}
Finally, we extended our analysis to the full, challenging 3D mask problem. The mask structure and parameters of the incident field are as in previous case for the WGNO. The nonuniform permittivities of the absorber in layer $j$ are described by:
\begin{align}
	&\eps_j(x,y){=} \frac{1}{4}\left[\tanh\left(\frac{x {+} a}{d}\right) {-} \tanh\left(\frac{x {-} a}{d}\right)\right] \left[\tanh\left(\frac{y {+} b}{d}\right) {-} \tanh\left(\frac{y {-} b}{d}\right)\right] (1 {-} \eps) {+} \eps,
\end{align}
where $\eps$ is permittivity of absorber of $j$th layer, $a=L_x / 4$, $b=L_y / 4$, $d = \lambda / 10$.

\begin{figure}[t!]\centering
	\includegraphics[width=0.5\textwidth]{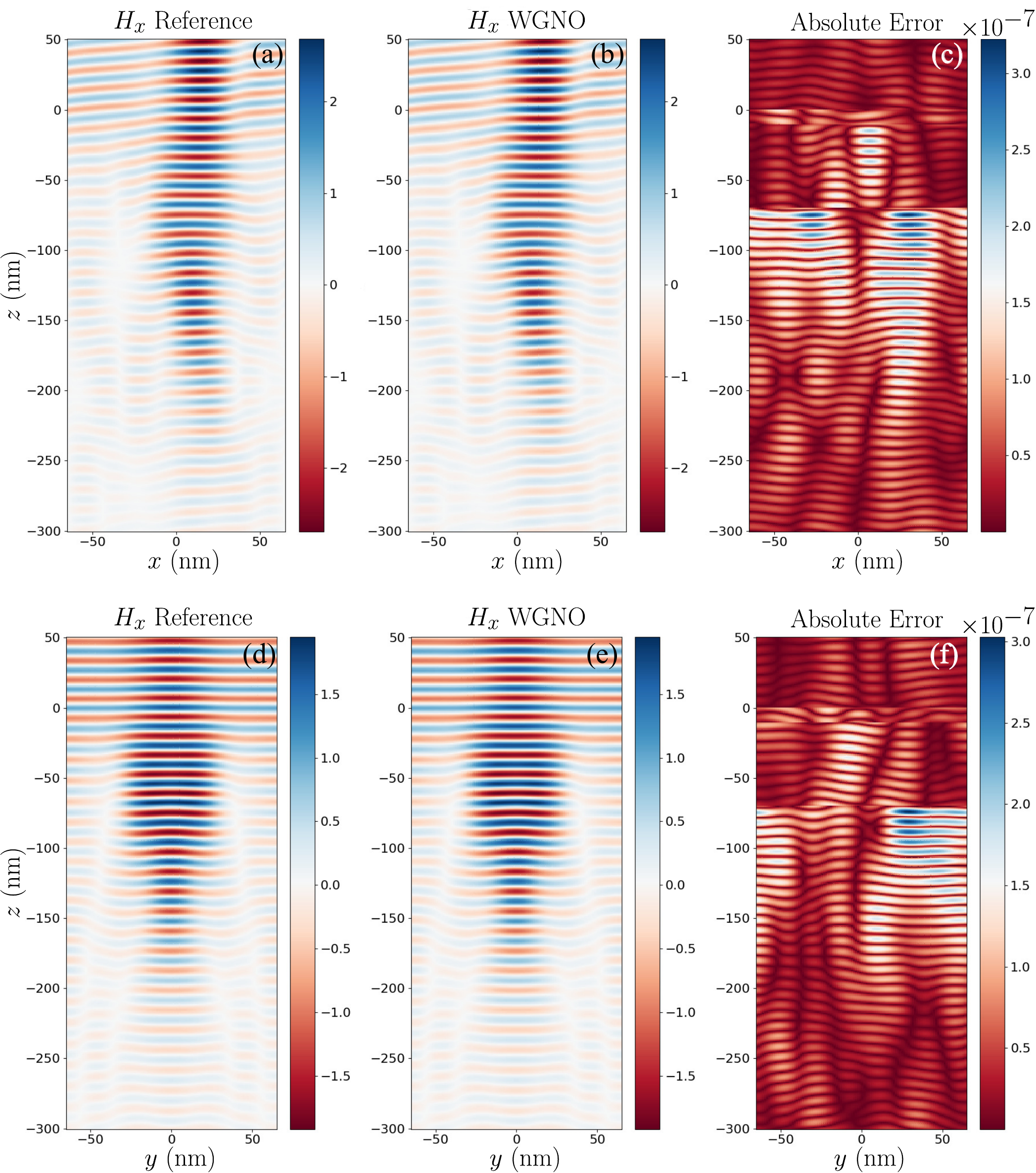}
	\caption{Comparison of reference WG solution ((a),(d)) and WGNO solution ((b),(e)) for a 3D mask at 13.5 nm, (c) and (f) are absolute errors. (a)--(c) are in cross-section $y=0$, (d)--(f) are in the cross-section $x=0$.}\label{fig:3d_mask_results}
\end{figure}

Fig.~\ref{fig:3d_mask_results} shows a cross-section of the field calculated by our WGNO compared to the reference solution at 13.5 nm. The agreement is excellent, demonstrating the capability of our method to handle the full 3D problem with high fidelity.

The performance metrics for the 3D case, shown in Table~\ref{tab:3d_mask}, are even more compelling. The WGNO maintains high accuracy while achieving a speedup of over 200 times compared to the rigorous WG solver. The entire training process for the 3D case took only about 18 seconds, demonstrating the remarkable efficiency and scalability of our proposed operator.

\begin{table}[h!]
	\centering
	\caption{Performance of WGNO on 3D lithography mask simulation. Training time is $18$ s.}
	\label{tab:3d_mask}
	\begin{tabular}{l|cc}
		\toprule
		Wavelength & Rel. $L_2$ Error & Inference (s) \\
		\midrule
		13.5 nm & $1.0\times 10^{-7}$ & $2.18\times 10^{-4}$ \\
		11.2 nm & $1.9\times 10^{-6}$ & $2.08\times 10^{-4}$ \\
		\bottomrule
	\end{tabular}
\end{table}

\section{Conclusion}
In this work, we have shown that physics-informed neural networks and neural operators can achieve high accuracy for complex diffraction problems of EUV lithography simulation. We determined that while training times can vary, the inference time is extremely small. For relatively thin masks, PINNs can provide solutions accurate enough for initial calculations. Most significantly, we have established that our proposed Waveguide neural operator is a state-of-the-art solver for this application. It provides highly accurate solutions for a full 3D mask, surpassing modern approaches in both accuracy and inference time, all while requiring only a very short training period. This level of performance makes the WGNO a highly promising tool for accelerating the OPC design cycle in semiconductor manufacturing.

\end {document}